\documentclass[11pt]{amsart}
\usepackage{a4wide,amssymb,epsfig,righttag}

\newtheorem{majorthrm}{Theorem}

\newtheorem{lemma}[equation]{Lemma}
\newtheorem{prop}[equation]{Proposition}

\theoremstyle{definition}

\theoremstyle{remark}
\newtheorem*{remark}{Remark}
\newtheorem*{remarks}{Remarks}
\newtheorem*{philosophical}{Philosophical remark}

\def\RE{{\mathbb R}}
\def\CX{{\mathbb C}}
\def\DD{{\mathcal D}}

\def\GG{{\mathcal G}}
\def\LD{{\mathcal L}}

\def\Ga{\Gamma}
\def\bX{\overline{X}}
\def\bg{\overline{g}}
\def\bc{\overline{c}}
\let\pf=\proof
\let\endpf=\endproof
\let\phi=\varphi
\let\epsilon=\varepsilon

\DeclareMathOperator{\Coker}{Coker}

\DeclareMathOperator{\End}{End}
\DeclareMathOperator{\ind}{index}
\DeclareMathOperator{\Image}{Im}
\DeclareMathOperator{\Ker}{Ker}

\begin{document}

\title{Gluing theorems for ASD metrics}
\author{A.G.\ Kovalev  and M.A.\ Singer}
\address{Department of Mathematics and Statistics, University of
Edinburgh}
\email{agk@maths.ed.ac.uk, michael@maths.ed.ac.uk}
\thanks{Preprint version, 6th February 1998}
\maketitle

In 1991, Andreas Floer showed how methods of non-linear elliptic PDEs
might be used to obtain anti-self-dual (ASD) metrics on the connected
sum of $l\geq 3$ copies of $\overline{\CX P^2}$ \cite{floer-conformal}.
Much subsequent work on 
ASD metrics has depended upon gluing theorems but has relied more on
the twistor-based methods pioneered by Donaldson and Friedman
\cite{donaldson-friedman}. (A very
important exception is Taubes's `stable existence theorem' for ASD metrics
\cite{taubes-conformal}.)

The purpose of this note is to report on work that goes in the
direction of realizing the full potential of the techniques introduced by
Floer. Full proofs will not be given here; the details, along with
generalizations and examples, will appear elsewhere \cite{glue2}. The
present exposition is intended to make accessible to a
wider audience the important insights that Floer brought to this
problem by stating some relatively non-technical, but representative
results.

Acknowledgements: The authors are grateful to a number of mathematicians
with whom they have discussed this work during the summer of 1997: Robin
Graham, Dominic Joyce, Claude LeBrun, Mario Micallef, Henrik Pedersen,
Yat-Sun Poon and Rick Schoen.  The second author is also grateful to the
IHES for its hospitality during November 1997; his stay there was very
helpful during the final stages of this work.
The first author was supported by a William Gordon Seggie Brown
Fellowship and the second author was supported by an EPSRC Advanced
Fellowship.

\section{ASD conformal structures}

Recall that a Riemannian metric $g$ on an oriented $4$-manifold $X$ is
said to be ASD if the self-dual part $W^+(g)$ of its Weyl tensor
vanishes. The condition is conformally invariant because
$W^+(\hat{g}) = W^+(g)$ for any two metrics $g$ and $\hat{g} = e^fg$
in the same conformal class $c$ on $M$; for this we regard $W^+$ as a
section of the bundle $E^2 \subset \Lambda^2T^* \otimes \End(T)$ which
corresponds, by raising an index, to $S^2_0\Lambda^+$. (For more
details and background, see for example
\cite{pen,ahs,taubes-conformal,taubes-glueing} or
\cite{gau} for the state of the art about 1992.)

The condition $W^+(g)=0$ is a non-linear PDE of second order in $g$
that is elliptic modulo the action of the diffeomorphism group of $X$.
More precisely, given an ASD metric $g$, there is an elliptic complex
(to be called the {\em deformation complex} \cite{king-kotschick})
\begin{equation} \label{dc}
C^\infty (X,E^0) \stackrel{L_g}{\longrightarrow}
C^\infty (X,E^1) \stackrel{D_g}{\longrightarrow}
C^\infty (X,E^2)
\end{equation}
where $E^0 = TX$, $E^1$ is the bundle of symmetric trace-free
endomorphisms of $TX$, and $E^2$ was defined in the previous
paragraph. A section of $E^0$ yields in the standard way an
infinitesimal diffeomorphism of $X$, while a section of $E^1$ yields a
tangent to the space of conformal structures on $X$ at $g$ via
$$
g(1 +th)(\xi,\eta) := g(\xi,\eta) + tg(\xi, h\eta)
$$
for small $t$. The operator $L_g$ gives the infinitesimal action of
the diffeomorphism group on the space of conformal structures---that
is, the Lie derivative, while $D_g$ gives the infinitesimal change in
$W^+$:
$$
L_g\xi = (\LD_\xi g)_0,\;\;\;
D_g h = \frac{d}{dt}W^+[g(1 +th)]|_{t=0}.
$$
($(\cdot)_0=$ trace-free part.)
With this choice of bundles, $L$ and $D$ are exactly conformally
invariant:
$$
L_{\hat{g}} = L_g,\;\;\;\; D_{\hat{g}} = D_g
$$
where $g$ and $\hat{g}$ are related as above.

   The cohomology groups of~\eqref{dc} will be denoted by $H^*_c(X)$
or just $H^*_c$ ($c$ standing for the conformal class of $g$), and
play a crucial role in the theory. In particular the vanishing of the
{\em obstruction space} $H^2_c$ is the basic hypothesis that enters
the statement of our gluing theorems.  Note that when $X$ is compact, the
$H^*_c$ are finite-dimensional vector spaces.

The entire discussion goes through in a standard way to include
oriented $4$-dimensional orbifolds, with smooth ASD orbifold
metrics. The significance of this is that such ASD orbifolds are
plentiful, cf.\ Remark~(ii) below. From now on we shall allow
ourselves to consider $4$-orbifolds, but only those with isolated
(i.e.\ codimension $4$) singular points.

Our first extension of Floer's work is the following:
	\begin{majorthrm}
	\label{thma}
For $i=1,2$, let $\bX_i$ be a  compact $4$-orbifold, with ASD
conformal structure $c_i$ and $H^2_{c_i}=0$. Suppose that $x_i$ is a
point of $\bX_i$ such that $x_1$ and $x_2$ are complementary.  Then 
the generalized connected sum $(\bX_1,x_1)\sharp (\bX_2,x_2)$ admits
ASD (orbifold) metrics.
\end{majorthrm}

As a point of notation, from now on we shall denote compact spaces
with overlines to distinguish them from the non-compact manifolds
which we shall be working with in \S\S2--3.  Moreover the
qualification `generalized' of the term connected sum will usually be
omitted.

To explain the word `complementary' in Theorem~\ref{thma} recall that
by definition of orbifold, there exists a neighbourhood $U_i$ of $x_i$
in $\bX_i$ and homeomorphisms $\phi_i: \RE^4/\Ga_i \to U_i$ (with
$\phi_i(0) = x_i$) where $\Ga_i \subset SO(4)$ is a finite subgroup,
such that the induced action of $\Ga_i$ on $S^3\subset(\RE^4 - 0)$ has
no fixed points. Then $x_1$ and $x_2$ are called complementary if
there exists an orientation-{\em reversing} element of $O(4)$ that
conjugates $\Ga_1$ into $\Ga_2$.  On this definition, any pair of
smooth points will be complementary as will any pair of 
singularities modelled on $\RE^4/\{\pm 1\}$.

The connected sum $(\bX_1,x_1)\sharp (\bX_2,x_2)$ is obtained by
deleting small balls centred on $x_i$ and joining the two boundaries
by a cylinder (or {\em neck}) $Y_1 \times [-l,l]$. Here $Y_i$ is the
link of $x_i$ in 
$\bX_i$; it will be a spherical space-form $S^3/\Ga_i$. It is clear that
for this to make sense, that is, for $(\bX_1,x_1)\sharp (\bX_2,x_2)$
to have an orientation compatible with the given orientations of the
$\bX_i$, one end of the cylinder must be joined by an
orientation-reversing map. The condition that the $x_i$ be
complementary precisely guarantees the existence of an
orientation-reversing isometry $Y_1 \to Y_2$. 
	\begin{remarks}
 (i) For ASD $4$-manifolds, this result was proved in
\cite{donaldson-friedman}, by twistor
methods. These methods were generalized to deal with orbifolds with
$\Ga_i=\{\pm 1\}$ in \cite{lebrun-singer}, and to more general cyclic
singularities by Jian Zhou \cite{jz}. 
A drawback
of the twistor approach is its increasing complication with the
complexity of the singularities of~$X_i$.  By contrast, the analytic
approach is insensitive to the complication of the
singularities---indeed the same methods prove Theorem~\ref{thmb} below,
in which, roughly speaking, the cross-section $Y$ of the neck is
allowed to be a rather general compact, oriented $3$-manifold.

(ii) The hypothesis $H^2_c=0$ is satisfied by many important classes
of ASD orbifolds: conformal compactifications of ALE gravitons,
$4$-dimensional quaternion-K\"ahler spaces, and ALE ASD--K\"ahler
spaces (such as the family constructed in \cite{counter}). The
vanishing proof for this last family of examples is a variant of the
analysis in \cite{exdef} and is based on unpublished work of LeBrun
and the second author; the details will appear elsewhere.

(iii) It follows from the proof that the metrics constructed by
Theorem~\ref{thma} themselves have vanishing obstruction spaces. That
is, Theorem~\ref{thma} remains true for connected sums of any
(finite) number of ASD orbifolds performed simultaneously along pairs
of complementary singularities.

(iv) If $H^2_{c_i}\not=0$, then the same methods go through to solve
`the infinite-dimensional part' of the ASD equations for metrics on
$\bX_1\sharp \bX_2$, and it may yet be possible to obtain ASD metrics
after some additional work \cite[\S6.4]{donaldson-friedman}.  In fact,
these same methods yield 
satisfactory answers to the natural questions about the moduli spaces
of ASD metrics on $(\bX_1,x_1)\sharp (\bX_2,x_2)$ that are raised by
the construction (cf.\ for example the discussion in~\cite[Ch.\ 7]{DK}
of the same issue in the case of ASD Yang--Mills connections).  
These points will be treated carefully in \cite{glue2}.  
\end{remarks}
The proof of Theorem~\ref{thma} involves the construction of an
approximately ASD metric on the connected sum and then an application
of the implicit function theorem (IFT) to find a nearby exactly ASD
metric.  The nub of the matter is appropriate choices of the geometric
construction of the connected sum and function spaces in which to
apply the IFT.  The crucial thing is to have the linear problem
(coming from (\ref{dc})) under uniform control as the diameter of the
neck of the connected sum shrinks to $0$. The issues involved are well
discussed in references such
as~\cite{DK,floer-conformal,taubes-conformal,taubes-glueing}. 
We shall therefore concentrate on obtaining control of the linear
problem, following Floer's elegant approach. This involves a
cylindrical model of the connected sum and the theory of elliptic
operators on manifolds with cylindrical ends developed by Lockhart and
McOwen.  This emphasis on carrying out all the analysis on non-compact
manifolds is in contrast to the approach employed by Donaldson
\cite{DK}, or Taubes \cite{taubes-conformal, taubes-glueing}, who work
as far as possible with {\em compact} spaces.

\section{Manifolds with cylindrical ends}

By a manifold with a cylindrical end (or CE-manifold) $X = X^-\cup_Y
X^+$, we mean a non-compact manifold $X$ decomposed as a compact
manifold $X^-$ with boundary $Y$, and a (half)-cylinder $X^+ = Y
\times [0,\infty)$, with $X_{\pm}$ attached along the common boundary
$Y=\partial X^-\cong Y\times\{0\}\subset X^+$. A CE-metric $g$ on $X$ is
one which approaches a Riemannian product metric $g_0$ on $X^+$ at an
exponential rate:
\begin{equation} \label{er}
\sup_{Y\times \{t\}}|g - g_0|_{g_{0}}  \leq Ce^{-\eta t},\;\;\;
\sup_{Y\times\{t\}}|\nabla^k g|_{g_0} \leq C_ke^{-\eta t},\;\;
\mbox{ for }t>0\mbox{ and }k=1,2,\ldots.
\end{equation}
Here the point-wise norms and the covariant derivative are those of
$g_0$ and $\eta>0$ is some constant.  The discussion extends in an
obvious way to deal with CE-orbifolds and CE-manifolds with more than
one end.  It will always be assumed, however, that $Y$ is smooth.

Given a compact orbifold $\bX$ with a marked point $x$ and smooth
orbifold metric $\bg$, one obtains a CE-manifold $X$ with CE-metric
$g$ by a conformal rescaling that is singular at $x$. Since the
transformation is conformal, $g$ is ASD whenever $\bg$ is so.  Indeed,
choose a small geodesic ball $B(r_0)$ about $x$, let $Y= \partial B=
S^3/\Ga$ be the 
link of $x$ in $\bX$ and let $r$ be the geodesic distance from $x$.
Then in $B$ the metric takes the form
$$
\bg = dr^2 + r^2(h_0 + r^2 h_2)
$$
where $h_0$ is the induced metric of constant positive curvature on $Y$
and $h_2$ is some smooth $r$-dependent family of metrics on $Y$. On $B
- x$ set
$$
g = r^{-2}\bg = dt^2 + h_0 + e^{-2t}h_2, \mbox{ where } r/r_0 = e^{-t}.
$$
In particular $g$ is a CE-metric with $\eta=2$ and $g_0$ the
standard product metric on $(S^3/\Gamma)\times [0,\infty)$.
To complete the formal description, $X^- = \bX
- B$.  In what follows,  we shall refer to $X$ as the
{\em conformal cylindrification} of $(\bX,x)$. (It is understood here that
$r$ is to be continued smoothly to a positive function on $X^-$.)

For another class of examples, for which we are indebted to Claude
LeBrun, let $Y$ be an oriented $3$-manifold of constant sectional
curvature and admitting an orientation-reversing isometric involution
$\iota$ with only isolated fixed points. Let $\hat{X} = Y \times \RE$
with the product metric; because $Y$ has constant curvature,
$\hat{X}$ is conformally flat. Now let $X = \hat{X}/\langle (\iota,-1)\rangle$
with the induced (orbifold) metric. Then $X$ is a conformally flat
CE-orbifold with 
just one end $X^+ = Y\times[0,\infty)$ and singular points
corresponding to the fixed points of $\iota$ on $Y\times\{0\}$.

\vspace{10pt}
Lockhart and McOwen \cite{lockhart-mcowen} gave a package of Fredholm
theorems for a class of elliptic operators on CE-manifolds; these are
the operators that are naturally adapted to the product geometry of
the cylindrical end, and will be called CE-operators\footnote{About
the same time Melrose and Mendoza \cite[Theorem~(51)]{Me} obtained
similar results through the development of a calculus of `totally
characteristic pseudo-differential operators'.  This appears to be
the more fruitful viewpoint in that there are important
generalizations to geometric situations in which the CE-model is
inappropriate \cite{Ma,Me}. For the purposes of this paper, we shall
stay with the CE framework.  Note however that
Proposition~\ref{excision} uses CE pseudo-differential operators and
hence relies on the Melrose--Mendoza theory.}.

By a CE (differential) operator $A:C^\infty(E) \to C^\infty(F)$ (where
$E$ and $F$ are smooth vector bundles over a CE-manifold~$X$) we mean
a (differential) operator which is asymptotic, at an exponential rate,
to a `time-independent' operator $B$ (say) on $Y \times \RE$. This
statement is understood relative to a choice of isomorphism of $E$
with the pull-back by the projection $Y\times [0,\infty) \to Y$ of
$E|Y \times 0$, similarly for $F$. Any differential operator
canonically defined by a CE-metric will automatically be a CE
differential operator; in particular, $L_g$ and $D_g$ in \eqref{dc}
have this property.

For the analysis of elliptic CE-operators, a crucial invariant is the
asymptotic spectrum $\Sigma(A)= \Sigma(B)$. This is the set of complex
numbers $\lambda$ such that
$$
B(e^{i\lambda t}u(y))= 0 \mbox{ has a solution for some }0\not=u \in
C^\infty (Y,E).
$$
It is a standard fact that $\Sigma(A)$ is a discrete subset of $\CX$
which meets every horizontal strip
$\delta_1 < \Image \lambda < \delta_2$
in a finite set  of points.  The other essential ingredient is the
introduction of the weighted Sobolev spaces $L^p_{k,\delta}$ defined
by completing the space of 
functions (or sections) with compact support in the
norm
$$
\|u \|_{L^p_{k,\delta}} = \sum_{r=0}^k \|e^{\delta t}\nabla^r u\|_p.
$$
Here the conventions are that $t$ is equal to the standard coordinate
on $X^+$ (as above) and is smoothly cut off to zero on the
compact piece $X^-$.
Second, $\nabla$ denotes a CE-covariant-derivative operator that
preserves a CE-bundle metric. Third the $L^p$-norm on the RHS is
calculated with this CE-metric on the bundle and a CE-metric on the
base $X$.

The basic result of \cite{lockhart-mcowen} is that an elliptic
CE-operator $A$ of order $m$ (say) extends to a bounded Fredholm map
$$
A_\delta : L^p_{k,\delta}(X,E) \to 
            L^p_{k-m,\delta}(X,F)
$$
if and only if $\delta$ is not the imaginary part of any $\lambda \in
\Sigma(A)$; accordingly we shall call $\delta$ an {\em
exceptional weight} (for $A$) if and only if $A_\delta$ is {\em not}
Fredholm. It will be convenient to denote by $\delta_0 = \delta_0(A)$
the first positive exceptional weight  for $A$.
As in the compact case, the index is independent of $p$
and $k$ and $\Ker(A_\delta)$ consists of smooth sections.  There is a
weighted version of the Fredholm alternative, identifying
$\Coker(A_\delta)$ with $\Ker(A^*_{-\delta})$ where $A^*$ is the
formal $L^2$-adjoint of $A$. The index depends
strongly upon $\delta$, however, and jumps according to the formula
\begin{equation} \label{jump}
\ind(A_\delta) - \ind(A_{\delta'}) = n(\delta,\delta') =
\sum_{\delta < \Image \lambda < \delta'} d(\lambda),
\end{equation}
for any non-exceptional $\delta < \delta'$
where $d(\lambda)$ is the dimension of the space of all elements of
the kernel of $B$ of the form 
\begin{equation} \label{kerfrm}
\exp(i\lambda t)\sum_0^N u_n(y)t^n.
\end{equation}
Notation: it will be convenient later also to write $\Ker_\delta(A) =
\Ker(A_\delta)$ and similarly for $\Coker$ and $\ind$.

\vspace{10pt}
We are now ready to formulate a generalization of Theorem~\ref{thma}.
Let $X$ be an ASD CE-manifold (or orbifold) with CE-metric $g$
satisfying \eqref{er}. Replace the complex (\ref{dc}) by the 
operator\footnote{This operator is elliptic of mixed order, but such
operators are treated explicitly in \cite{lockhart-mcowen}.}
$$
\DD_g = (D_g,L^*_g): C^\infty(X,E^1) \to C^\infty(X,E^2) \oplus C^\infty(X,E^0)
$$
($L^*_g$ being the formal adjoint of $L_g$). On a compact manifold the
kernel and cokernel of $\DD_g$ would of course be respectively
isomorphic to $H^1_c$ and $H^2_c \oplus H^0_c$; in the CE case, pick
any weight $\delta>0$ less than $\min(\eta,\delta_0(\DD_g))$ and set
\begin{equation} \label{dcd}
H^1_g = \Ker_\delta(\DD_g),\;\;H^2_g\oplus H^0_g =
\Coker_\delta(\DD_g).
\end{equation}
We may now state
	\begin{majorthrm}
	\label{thmb}
For $i=1,2$, let $X_i = X_i^- \cup_{Y_i} X_i^+$ be ASD CE-orbifolds
with CE-metrics $g_i$ and suppose there exists an
orientation-reversing isometry $\iota:Y_1 \to Y_2$. Then if
$H^2_{g_i}=0$,  there exist ASD metrics on the glued orbifold
	\begin{equation}
	\label{union}
X(l) = X_1^-\cup_{Y_1} Y_1 \times [-l,l] \cup_{Y_2} X_2^-
	\end{equation}
for all $l$ sufficiently large. 
	\end{majorthrm}
Here
 the role of $\iota$ is to glue the right-hand end of the cylinder
$Y\times [-l,l]$ onto the boundary of $X^-_2$ so that $X(l)$ has an
orientation compatible with the given orientations of the $X_i$ (cf.\
the discussion after Theorem~\ref{thma}).   We
have not placed any geometric conditions on the $Y_i$ in this
statement. However, the decay condition~\eqref{er} clearly entails the
existence 
of a metric $h$ on $Y$ such that the product metric on $Y\times \RE$
is ASD, and one may show that any such ASD product metric is
conformally flat. This in turn forces $h$ to be a metric of constant
sectional curvature.  
	\begin{philosophical} 
One might have expected an additional hypothesis concerning the
 vanishing of `obstructions from the
neck' to enter in the statement of 
this theorem.  Although this does not occur
explicitly, note that in general $ H^2_g = \Ker_{-\delta}(D^*_g) $ may contain
elements that do not decay along the cylinder $Y \times
[0,\infty)$.  In this sense, it seems appropriate to think of the
hypothesis $H^2_g =0$ as {\em including} the hypothesis that there
are no obstructions from the neck.  We make no attempt here to make
 this notion more precise, however.

	\end{philosophical}

\begin{remark} Remarks (iii) and (iv) following Theorem~\ref{thma}
apply, {\em mutatis mutandis}, to this Theorem also. (In particular,
$X_i$ and $Y_i$ are allowed to have multiple components.)
\end{remark}

Theorem~\ref{thma} is contained in Theorem~\ref{thmb} by taking $X_i$
to be the `conformal cylindrification' of $(\bX_i,x_i)$. Then the
existence of $\iota$ is equivalent to the complementarity of $x_1$ and $x_2$.
The remaining
point is to compare the cohomology groups of (\ref{dc}) with
(\ref{dcd}).  This is achieved by the following 
	\begin{majorthrm}
  Let $(\bX,x)$ be a compact ASD $4$-orbifold with metric $\bg$ and
let $X$ be the conformal cylindrification, with metric 
$g$. Then we have
$$
H^*_g(X) = H^*_{\bc}(\bX,x)
$$
where on the RHS we have the cohomology of the complex
\begin{equation} \label{dcf}
C^\infty_x(\bX,E^0) \stackrel{L_{\bg}}{\longrightarrow}
C^\infty(\bX,E^1) \stackrel{D_{\bg}}{\longrightarrow}
C^\infty(\bX,E^2)
\end{equation}
and the subscript denotes sections which vanish at $x$.
\label{thmc}\end{majorthrm}
Since  $H^2_{\bc}(\bX,x)= H^2_{\bc}(\bX)$, it follows at once that the
hypotheses regarding the vanishing of $H^2$ in Theorems~\ref{thma}
and \ref{thmb} agree.
	\begin{remark}
 Theorem~\ref{thmb} may be viewed as the natural basic
gluing theorem for ASD metrics in `the CE category'.  We have seen
that any compact ASD orbifold $\bX$ with a marked point $x$ can be
regarded as an ASD CE orbifold.  With this in mind, our main results
can be summarized by:
\begin{quote}
Theorem~\ref{thmb} $+$ conformal cylindrification $+$
Theorem~\ref{thmc} $\Longrightarrow$ Theorem~\ref{thma}.
\end{quote}
This organization of the material seems to be the logical conclusion
of Floer's methods and exposes most clearly the
structure of the argument. In any case,
Theorem~\ref{thmb} and its proof are of independent interest.
\end{remark}

\section{Structure of the proofs of Theorems \ref{thma} and \ref{thmb}}

Theorem~\ref{thmb} and, accordingly, Theorem~\ref{thma} produce
solutions of a non-linear elliptic PDE on a manifold~$X$ from those on
its component pieces~$X_i$, $i = 1,2$. The arguments therefore involve
applications of analysis---carefully tailored to the geometry
of~$X$. More specifically, the theme is the comparison of the
linearized problem on $X(l)$ and on its components, corresponding to
decomposition~\eqref{union}.  By way of preparation, before explaining
the main ingredients in the proof we need to introduce some notation.
\medskip

Suppose, for simplicity, that in the statement of Theorem~\ref{thmb}
the $X_i^-$ are smooth compact manifolds with boundary, and the $X_i$
have just one cylindrical end.  First we construct an approximately
ASD metric $g(l)$ on $X(l)$. For this, fix a  standard
cut-off function $\alpha(t):\RE \to [0,1]$, equal to $1$ for $t\leq
0$ and $0$ for $t\geq 1$. Denote by $t_i$ the standard coordinate
along the cylinder $X_i^+$. By cutting off the exponentially decaying
part of $g_i$ with $\alpha(t_i-l+1)$ we obtain a smooth metric
$g_i(l)$ on $X_i$ which is equal to the cylinder metric for $t_i \geq
l$. Therefore the map
\begin{equation} \label{attach}
(y,t_1)
 \longmapsto (\iota(y), 2l-t_2)
\end{equation}
is an orientation-preserving isometry from a small neighbourhood of
$\{t_1 =l\}$ to a corresponding neighbourhood of $\{t_2 =l\}$. Let
$X(l)$ be the compact manifold obtained by attaching $X_1$ to $X_2$ by
\eqref{attach}, and denote by $g(l)$ the obvious metric induced from
the $g_i(l)$ on $X(l)$. 
\begin{figure}[hbt]
\centering
\epsfig{file=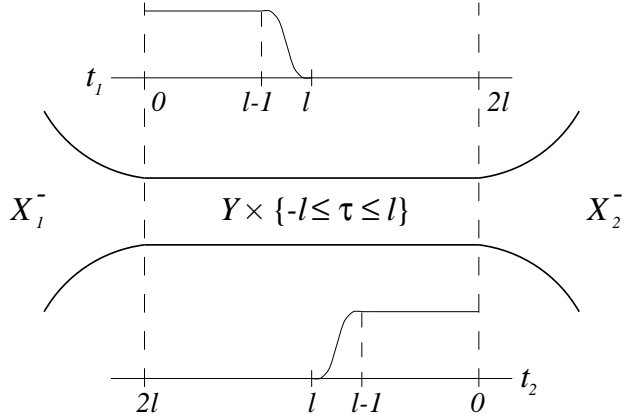}
\caption{The manifold $X(l)$ with the cut-off
functions used to construct $g(l)$.}\label{the-picture}
\end{figure}
In addition to the $t_i$ it is useful to introduce the coordinate
$\tau$ on the neck 
$Y_1 \times [-l,l]$
of $X(l)$, so that $\tau = t_1-l = l -t_2$. Then,
for each large enough $l$, we have $\tau=0$ in the middle of the neck,
and the neck itself is given as the region $|\tau| \leq l$ and thus
has length $2l$ (Fig.~\ref{the-picture}). 

It is easy to see that for each $p$ and $k$, we have
	\begin{equation}
	\label{firstest}
\| W^+[g(l)]\|_{p,k} \leq C_{p,k} \exp(-\eta l)
	\end{equation}
where the $L^p_k$-norm on $X(l)$ is defined by $g(l)$. In other
words, the $g(l)$ form an improving sequence of approximately ASD
metrics on $X(l)$, as $l \to \infty$.  To fit in with the CE-Fredholm
theory outlined in \S 2, however, we need a variant of
\eqref{firstest}, involving the introduction of a weight-function
$w(l)$ on $X(l)$. 
Put $w(l)$ equal to unity on $X^-_1$ and $X^-_2$, and extend it to a
smoothed version of the function $\exp(\delta(l-|\tau|))$ on the neck.
Define $L^p_{k,w(l)}(X(l))$ to be the $w(l)$-weighted Sobolev space
using the metric $g(l)$.  Of course, since $X(l)$ is compact, for any
fixed $l$ this
weighted norm is equivalent to the
unweighted $L^p_k$-norm; but since $w(l) \approx e^{\delta l}$
near $\tau=0$, this equivalence is {\em not} uniform as $l\to
\infty$.  Let
	\begin{equation}
	\label{three-pieces}
U(l) = L^p_{2,w(l)}(X(l),E^1),\;\;\;
V(l) = L^p_{0,w(l)}(X(l),E^2),\;\;\;
W(l) = L^p_{1,w(l)}(X(l),E^0),
	\end{equation}
and denote by 
$$
\DD(l): U(l)  \longrightarrow V(l) \oplus W(l)
$$
the extension of $(D_{g(l)}, L^*_{g(l)})$ to a bounded operator between
these spaces.  Then~\eqref{firstest} implies 
	\begin{equation}
	\label{secest}
\| W^+[g(l)]\|_{V(l)} \leq C \exp(-(\eta-\delta)l),
	\end{equation}
so that if $\delta < \eta$ the $g(l)$ still form an improving sequence of
approximately ASD metrics, as measured by the norm of $V(l)$.

Returning to~\eqref{three-pieces}, consider three operators $\DD_0$,
$\DD_1$ and $\DD_2$, 
$$
\DD_a: U_a \to V_a\oplus W_a, \qquad   a=0,1,2,
$$
where
$$
U_i = L^p_{2,\delta}(X_i,E^1),\;\;\;
V_i = L^p_{0,\delta}(X_i,E^2),\;\;\;
W_i = L^p_{1,\delta}(X_i,E^0),
$$
and
$$
U_0 = L^p_{2,w}(X_0,E^1),\;\;\;
V_0 = L^p_{0,w}(X_0,E^2),\;\;\;
W_0 = L^p_{1,w}(X_0,E^0),
$$
where $X_0 = Y \times \RE$ (with coordinate $\tau$ on $\RE$) and $w$
is (a smoothed version of) the function $e^{-\delta|\tau|}$ on $X_0$.
The aim is now to understand $\DD(l)$ in terms of the  $\DD_a$.

Write
$$
H_a = \Ker(\DD_a),\;\;\;J_a\oplus K_a = \Coker(\DD_a)
$$
(relative to the direct-sum decomposition $V_a\oplus W_a$).
Now
$J_0\oplus K_0=0$ since this cokernel can be identified with
$\Ker_{w^{-1}}(\DD_0^*)$; and it is impossible for a non-zero element
of the kernel of a $\tau$-independent operator to decay exponentially
at both ends of $Y\times \RE$ (this is an easy consequence of the
results of~\cite{kondratiev,maz'ya-plamenevskij}).  
In terms of previous notation,
$$
H_i = H^1_{g_i},\;\;\;J_i = H^2_{g_i},\;\;\;K_i = H^0_{g_i}.
$$
Thus the hypothesis of Theorem~\ref{thmb} is $J_i=0$; for simplicity,
we shall assume also that $K_i=0$. (One can always reduce to the case
$K_i=0$ by the device of framing the problem at a finite number of
points, precisely as one can deal with the reducible connections in gauge
theory.) 
\medskip

The comparison of $\DD(l)$ with the $\DD_a$ involves a topological
idea and an analytical idea. The topology enters through index theory;
in \S\ref{excis}, an appropriate version of the excision property of
the index is used to obtain the simple formula
	\begin{equation}
	\label{if}
\ind(\DD(l))  = \ind(\DD_0) + \ind(\DD_1) + \ind(\DD_2).
	\end{equation}
The analytical step  involves the construction of a subspace
$U^\perp(l)\subset U(l)$
with the following properties:
	\begin{enumerate} \renewcommand{\labelenumi}{(\roman{enumi})}
	\item
The {\em main estimate} holds: there exists $C>0$ such that for
all $l \geq l_0$, and $h \in U^\perp(l)$,
	\begin{equation}
	\label{me}
\|\DD(l)h\|_l  \geq C\|h\|_l;
	\end{equation}
	\item
$U^\perp(l)$ is of the `correct' codimension:
	\begin{equation}
	\label{cc}
U(l)/U^\perp(l) \cong H_0 \oplus H_1 \oplus H_2.
	\end{equation}
	\end{enumerate}

(In the statement of~\eqref{me} we have used an obvious notational
simplification for the norms involved.)

It now follows from \eqref{if}, \eqref{cc} and \eqref{me} that the
restriction of $\DD(l)$ to $U^\perp(l)$ is surjective, with
uniformly bounded inverse
$\GG(l)$, say, satisfying
$\|\GG(l)\| \leq C^{-1}$ for $l \geq l_0$.  This, together with
\eqref{secest}, is precisely what is needed for a successful
application of a standard modification of the IFT (cf.\ \cite[Lemma
(7.2.23)]{DK} or \cite[Lemma 4.2]{floer-conformal}).

\subsection{The excision property} \label{excis}

The main idea behind the proof of \eqref{if} is the following claim.
	\begin{prop}[excision property]
	\label{excision}
Let $X$ be a CE-manifold, $E$ and $F$ vector bundles over
$X$ and $A:C^\infty(E) \to C^\infty(F)$ an elliptic differential
CE-operator. Suppose that there exists an open subset $U$ of $X$ and
trivializations
$$
\alpha : E|U \to U \times \CX^n,\;\;\;
\beta : F|U \to U \times \CX^n.
$$
Then for every $V \Subset U$ and $D>0$ there exists an elliptic CE
pseudo-differential operator \mbox{$P:C^\infty_0 (E) \to C^\infty(F)$}
(of order zero) such that
	\begin{enumerate}
	\renewcommand{\labelenumi}{(\alph{enumi})}
	\item
for $|\delta| < D$, $A_\delta$ is Fredholm if and only if $P_\delta$ is
Fredholm, and $\ind(P_\delta) = \ind(A_\delta)$;
	\item
$P$ is equal to the identity over $V$, i.e. for every $u \in
C^\infty(E)$ with support contained in~$V$, $Pu =
\beta^{-1}\alpha\,u$.
	\end{enumerate}
	\end{prop}
This is an extension of
the usual excision principle for elliptic operators~\cite[\S8]{as}
to the non-compact setting of the CE-category. 
In order to motivate its use in deriving \eqref{if} let us consider
how the standard excision principle yields a relative index formula
for the change in the index under the operation of connected sum.

To be specific,
let $I(\bX)$ denote the index of \eqref{dc} over a compact 4-manifold
$\bX$,
$$
I(\bX) := \dim H^1_{\bc} - \dim H^0_{\bc} -\dim H^2_{\bc}  .
$$
Given compact smooth 4-manifolds, $\bX_1, \bX_2, \bX_1', \bX_2'$ say,
the standard excision principle yields a formula
	\begin{equation}
	\label{grif}
I(\bX_1\sharp \bX_2) - I (\bX_1\sharp \bX_2') 
= I( \bX_1'\sharp \bX_2) - I( \bX_1'\sharp \bX_2')  ,
	\end{equation}
since the bundles involved are trivial over the neck $V$ of each connected
sum. (Cf.~\cite[\S7.1]{DK} for an analogous application in gauge theory.)
Now take $\bX_1'=  \bX_2' = S^4$, to obtain
	\begin{equation}
	\label{rif}
I(\bX_1\sharp \bX_2) = I(\bX_1) + I(\bX_2) - I(S^4).
	\end{equation}
The formula~\eqref{rif} extends to the (generalized) connected sums of
orbifolds $(\bX_1,x_1)\sharp (\bX_2,x_2)$, with $S^4$ replaced by
$S^4/\Ga$, where $\Ga$ is the local isotropy of $x_1$ (or $x_2$).

Now consider the situation of interest in Theorem~\ref{thmb} and
\eqref{if} where one has a `generalized connected sum', the neck
having cross-section $Y$. Since the tangent bundle of any compact,
oriented $3$-manifold is trivial and the $E^i$ are all associated to
the tangent bundle of the ambient $4$-manifold, the restrictions of 
the $E^i$ to
subsets of the form  $Y\times I$ of $X(l)$ and of the $X_i$ are all
trivial.  One deduces from this observation and
Proposition~\ref{excision} a generalization 
of \eqref{grif} in which some or all of the compact manifolds
$\overline{X}$ may be replaced by CE-manifolds, $I(\overline{X})$
being replaced by $\ind_{\delta}(\DD_g)$ for some appropriate choice
of $\delta$.  In particular, replacing $S^4$ in~\eqref{rif} by the 
half-cylinder $Y\times [0,\infty)$, we obtain
$$
\ind (\DD(l)) = \ind (\DD_1) + \ind(\DD_2) - \ind_{w^{-1}}(\DD_0)
$$
where $w^{-1}$ is the weight $e^{\delta|t|}$ on $Y\times \RE$. To
prove \eqref{if}, it remains, therefore, to show the equality
$$
\ind_{w^{-1}}(\DD_0) = - \ind_{w}(\DD_0) .
$$
But, for any non-exceptional $\delta$, the operator $\DD_0$ is
invertible on the $\exp(\delta t)$-weighted Sobolev spaces over the
full cylinder~\cite{kondratiev,maz'ya-plamenevskij}; see also
\cite{lockhart-mcowen}. In particular $\ind_\delta(\DD_0) = 0$. Then,
\eqref{jump} gives 
$$
\ind_{w^{-1}}(\DD_0) - \ind_{\delta}(\DD_0) = - n(-\delta,\delta) = - 
[\ind_{w}(\DD_0) - \ind_{\delta}(\DD_0) ]
$$
and this completes the proof of \eqref{if}.

\subsection{The main estimate}

First we explain how to construct $U^\perp(l)$.  It will be defined
by 3 types of orthogonality conditions, corresponding to $H_0$,
$H_1$, $H_2$.  The obvious transversal to $H_i$ in $U_i$ is the $L^2$
orthogonal complement (with respect to $g_i$) of $H_i$.  For
technical reasons we modify this: since $H_i$ consists of
exponentially decaying sections, there exists $L>0$ such that the
orthogonal complement of $\tilde{H}_i$ is also transverse to $H_i$,
where $\tilde{H}_i$ consists of elements of the form
$\alpha(t_i-L)e_i$, where $e_i\in H_i$.   We know also that $H_0$ is
spanned by sections over $X_0$ of the form \eqref{kerfrm}, so we can
choose a space $\tilde{H_0}$ of the same dimension as $H_0$, but
consisting of sections supported in $\{|\tau| \leq \varepsilon\}$ and
such that $\tilde{H_0}^\perp$ is transverse to $H_0$. 
Transferring these conditions in the obvious way to $X(l)$ (for
$l\geq L+1$) we define $U^\perp(l)$ to consist of $h\in U(l)$ such that
	\begin{subequations}
	\label{ort}
	\begin{alignat}{2}
&h \bot e_i,  
&&\qquad\text{for any }e_i\in \widetilde{H}_i, \; i=1,2,
	\label{orta} \\
&h \bot \tilde{h}_0,  
&&\qquad\text{for any }\tilde{h}_0\in\tilde{H}_0,
	\label{ortb}
	\end{alignat}
	\end{subequations}
Since the supports of the $\alpha_i(l)$ are disjoint and do not meet
$\{|\tau| < \varepsilon\}$, 
it is clear that $U^\perp(l)$ satisfies~\eqref{cc}. 

Remark that in the basic example of $Y=S^3/\Ga$, one can show that
all elements of 
$H_0$ are independent of $\tau$. Then condition~\eqref{ortb} can be
defined more naturally by orthogonality conditions on the
restriction of $h$ to $Y \times \{\tau=0\}$ \cite[Eqn
(4.5) (1)]{floer-conformal}.

The main estimate is proved as follows. 
If~\eqref{me} fails, then there is a sequence $h_n \in U^\perp(l_n)$
with $l_n\to\infty$ and such that
	\begin{equation}
	\label{wrong}
\|h_n\|_n = 1,\;\;
\|\DD(l_n) h_n\|_n \to 0 \mbox{ as }n\to\infty.
	\end{equation}
(Here again we are using obvious notational simplifications.)  The
first step is to obtain control of the $h_n$ near the middle of the
neck, $\tau=0$.  More precisely, we have:
\begin{lemma} Given a sequence $h_n$ satisfying \eqref{wrong}, there
exists a subsequence $h_{n_j}$ with the following property. Given
$T > \varepsilon$ (as in the definition of
$\widetilde{H}_0$), let $K = Y \times \{|\tau| \leq T\}
\subset X(l)$ for  $l>T$; then
\begin{equation} \label{neckcontrol}
\lim_{j\to\infty} \|\exp(\delta l_{n_j}) h_{n_j}\|_{L^p_2(K,E^1)} =0.
\end{equation}
\end{lemma}

\pf Observe first that on $K$ as in the lemma, every coefficient of
$\DD(l_n)-\DD_0$ decays like $\exp(-\eta l_n)$ as $n\to\infty$. Note
similarly that the $L^p_{2,w(l_n)}$-norm of $h_n|K$ is uniformly comparable
to the $L^p_2(K)$-norm of $\exp(\delta l_n) h_n$ for all $n$. Combining
these two with the basic elliptic estimate for $\DD_0$, we obtain
$$
\|\exp(\delta l_n) h_n\|_{L^p_2(K)} \leq C(K)(\|\DD_n h_n\|_n + 
\|h_n\|_{L^p_{0,w(l)}(K)}),
$$
with $C(K)$ independent of $n$. So given \eqref{wrong}, the lemma will
be proved if we find a subsequence $h_{n_j}$ such that
\begin{equation} \label{wts}
\lim_{j\to\infty} \|\exp(\delta l_{n_j}) h_{n_j}\|_{L^p(K)} =0
\end{equation}
For this, define $h_n^{(0)} = \psi_n\exp(\delta
l_n)h_n$, where $\psi_n = \alpha(|\tau|- l_n)$. Regard $h_n^{(0)}$ as
a sequence of elements of $U_0$ by identifying the regions $|\tau|<
l_n$ of $X(l)$ and $X_0$.  From the
definitions of the weights $w(l)$ and $w$ and Eqn.~\eqref{wrong} this
is a bounded sequence in  $U_0$; hence there is a
subsequence such that $h_{n_j}^{(0)}$  converges weakly to
$h^{(0)}_\infty \in U_0$. With $T$ and $K$ as before, put
$f_j =
h^{(0)}_{n_j}|K$ and $f_\infty = h^{(0)}_\infty|K$. Using again the two
observations made at the beginning of the proof this weak convergence
is enough to give $\DD_0 f_\infty=0$. Since $T$ was
arbitrary, $\DD_0 h_\infty^{(0)}=0$. But the conditions
\eqref{orta} are preserved in the weak limit, so $h^{(0)}_\infty =0$.

Reformulating this slightly, we conclude that $\exp(\delta
l_{n_j})h_{n_j}$ is weakly convergent to $0$ in $L^p_2(K)$. Since $K$
is compact this gives strong convergence to $0$ in $L^p_1(K)$, proving
\eqref{wts}. \endpf

   An analogous, but simpler argument controls the
behaviour of $h_n$ over the $X_i$. With $h_n$ now denoting the
subsequence given by the lemma, define $h_n^{(i)} =
\alpha_{i,n}h_n$, ($\alpha_{i,n} = \alpha(t_i - l_n + 2)$)
 and regard $h_n^{(i)}$ as a sequence of sections over $X_i$.
By \eqref{wrong}, the sequence is bounded in $U_i$, and clearly
satisfies \eqref{ortb} once $l_n >L+1$.  We have
$$
\DD_i h^{(i)}_n = \alpha_{i,n}\DD_i h_n + [\DD_i,\alpha_{i,n}] h_n
   = \alpha_{i,n}\DD(l_n) h_n + [\DD_i,\alpha_{i,n}] h_n
$$
since $\DD(l_n) = \DD_i$ on the support of $\alpha_{i,n}$. As $n \to
\infty$, both terms on the RHS tend to zero (in the norm of $V_i
\oplus W_i$); the first because of \eqref{wrong}, the second
because of \eqref{neckcontrol}---for the operator
$[\DD_i,\alpha_{i,n}]$ is supported in the neighbourhood $Y\times
[-2,2]$ of the middle of the neck. Thus
\begin{equation} \label{bodycontrol}
\lim_{n\to \infty}\|h^{(i)}_n\|_{U_i} =0.
\end{equation}
Combining \eqref{neckcontrol} and \eqref{bodycontrol} contradicts the
first part of \eqref{wrong}, as required. 

\subsection{The proof of Theorem~\ref{thmc}}

Unlike the excision formula and the main estimate, which are very
general features of the behaviour of linear elliptic operators over
manifolds with long necks, Theorem~\ref{thmc} exploits specific
features of the ASD equations, above all their conformal invariance.
The argument is facilitated by making the special choice $\delta = 2 -
4/p$, $0 < \delta < 1$ (so $2 <p <4$).  Then, the natural inclusion
 $j:X \hookrightarrow \bX$ induces the following Banach-space isomorphisms:
\begin{equation} \label{cfw}
L^p(\bX, E^2) =L^p_\delta(X,E^2),\;\;\;
\{h\in L^p_2(\bX,E^1): h(x) =0\} = L^p_{2,\delta}(X,E^1).
\end{equation}
Since $j$ is merely a {\em conformal} map, it is
essential to keep track of the `conformal weights' of the bundles
$E^i$ and the way in which the point-wise and global norms change when
pulled back by $j$.  The first isomorphism of \eqref{cfw} may be
thought of as a generalization of the familiar fact that the
(global) $L^2$-norm of the $W(g)$ (or of $W^{\pm}(g)$) is a conformal
invariant; that is the case $\delta=0$, $p=2$.  By contrast $E^1$ has
conformal weight zero and its point-wise norm is conformally
invariant.  It follows that the two spaces in the second isomorphism
of \eqref{cfw} may be treated as Sobolev spaces of ordinary
functions, and then the equality is due to Biquard~\cite{Biquard}.

The isomorphisms in \eqref{cfw} make it straightforward to compare the
second cohomology groups.  The comparison of zeroth cohomology is
quite elementary and comparison of first cohomology follows
\cite[Prop.\ 3.2]{floer-conformal}, together with an application of
Proposition~\ref{excision}. Full details will appear in \cite{glue2}.

\subsection{Completion of proof of Theorem~\ref{thmb}}

As mentioned at the end of \S1, a version of the implicit function
theorem is used to go from the approximately ASD metric $g(l)$ to a
genuine ASD metric.  To be more precise, one uses the IFT to solve the
equations
\begin{equation} \label{nonlinear}
W^+[g(l)(1 + h(l))] = 0,  L^*_{g(l)} h(l) = 0
\end{equation}
such that the $U(l)$-norm of $h(l)$ is small.  In $4$ dimensions
$L^p_2 \subset C^0$ if $p>2$, so a small $U(l)$-norm ensures that
$\tilde{g}(l):=g(l)(1 +h(l))$ is a non-singular $L^p_2$-metric on
$X(l)$.  It then follows from an elliptic regularity argument with
\eqref{nonlinear} that $\tilde{g}(l)$ is smooth.  These remarks
show why it is necessary to take $p>2$ in our choice of function
spaces.


\end{document}